\documentclass[hidelinks,onefignum,onetabnum]{siamart220329}

\usepackage{graphicx}%
\usepackage{multirow}%
\usepackage{amsmath,amssymb,amsfonts,txfonts,newtxmath}%
\usepackage{mathrsfs}%
\usepackage[title]{appendix}%
\usepackage{xcolor}%
\usepackage{textcomp}%
\usepackage{manyfoot}%
\usepackage{booktabs}%

\usepackage{listings}%
\usepackage[capitalize]{cleveref}
\usepackage{braket}
\usepackage{nicefrac}

\usepackage{stackengine}
\usepackage{braket}
\usepackage{amsmath}
\usepackage{amscd}
\usepackage{amsfonts}
\usepackage{amssymb}
\usepackage{subfigure}
\usepackage{tabularray}
\UseTblrLibrary{booktabs}
\usepackage{hyperref} 
\usepackage{pgfplots}
\pgfplotsset{compat=newest}
\usetikzlibrary{calc}
\usepackage{tikzscale}
\usepackage{tikz-3dplot}
\usepackage[capitalize]{cleveref}
\crefname{equation}{}{}
\usepackage{nicefrac}

\usepackage[textsize=small,obeyFinal]{todonotes}

\usepackage{pgfplots}
\pgfplotsset{compat=newest}
\usetikzlibrary{calc, arrows.meta, bending, patterns}
\usepackage{tikzscale}
\usepackage{tikz-3dplot}

\newcommand{\logLogSlopeTriangle}[6]
{
	
	\pgfplotsextra
	{
		\pgfkeysgetvalue{/pgfplots/xmin}{\xmin}
		\pgfkeysgetvalue{/pgfplots/xmax}{\xmax}
		\pgfkeysgetvalue{/pgfplots/ymin}{\ymin}
		\pgfkeysgetvalue{/pgfplots/ymax}{\ymax}
		
		\pgfmathsetmacro{\xArel}{#1}
		\pgfmathsetmacro{\yArel}{#3}
		\pgfmathsetmacro{\xBrel}{#1-#2}
		\pgfmathsetmacro{\yBrel}{\yArel}
		\pgfmathsetmacro{\xCrel}{\xArel}
		
		\pgfmathsetmacro{\lnxB}{\xmin*(1-(#1-#2))+\xmax*(#1-#2)} 
		\pgfmathsetmacro{\lnxA}{\xmin*(1-#1)+\xmax*#1} 
		\pgfmathsetmacro{\lnyA}{\ymin*(1-#3)+\ymax*#3} 
		\pgfmathsetmacro{\lnyC}{\lnyA+#4*(\lnxA-\lnxB)}
		\pgfmathsetmacro{\yCrel}{(\lnyC-\ymin)/(\ymax-\ymin)} 
		
		\coordinate (A) at (rel axis cs:\xArel,\yArel);
		\coordinate (B) at (rel axis cs:\xBrel,\yBrel);
		\coordinate (C) at (rel axis cs:\xCrel,\yCrel);
		
		\draw[#5]   (A)-- 
		(B)-- 
		(C)-- node[pos=0.5,anchor=west] {$#6$}
		cycle;
	}
}

\usepackage{lipsum}
\usepackage{amsfonts}
\usepackage{graphicx}
\usepackage{epstopdf}
\usepackage{algorithmic}
\ifpdf
  \DeclareGraphicsExtensions{.eps,.pdf,.png,.jpg}
\else
  \DeclareGraphicsExtensions{.eps}
\fi

\newsiamremark{remark}{Remark}
\newsiamremark{assumption}{Assumption}
\newsiamremark{hypothesis}{Hypothesis}
\crefname{hypothesis}{Hypothesis}{Hypotheses}
\newsiamthm{claim}{Claim}

\headers{FEM for non-Lipschitz semilinear equation}{B. Vexler}

\title{A priori error estimates for finite element discretization\\ of semilinear elliptic equations\\ with non-Lipschitz nonlinearities}

\author{Boris Vexler\thanks{Technical University of Munich,
		School of Computation,
		Information and Technology,
		Department of Mathematics
  (\email{vexler@tum.de}, \url{https://www.math.cit.tum.de/math/vexler}).}
}

\usepackage{amsopn}

\ifpdf
\hypersetup{
  pdftitle={FEM for non-Lipschitz semilinear equation},
  pdfauthor={B. Vexler}
}
\fi

\begin{document}

\maketitle

\begin{abstract}
In this paper we develop numerical analysis for finite element discretization of semilinear elliptic equations with potentially non-Lipschitz nonlinearites. The nonlinearity is essecially assumed to be continuous and monotonically decreasing including the case of non-Lipschitz or even non-Hölder continuous nonlinearities. For a direct discretization (without any regularization) with linear finite elements we derive error estimates with respect to various norms and illustrate these results with a numerical example.  
\end{abstract}



\newcommand{\Lap}{\upDelta}
\newcommand{\R}{\varmathbb{R}}
\newcommand{\M}{{\mathcal M}}
\newcommand{\sgn}{\operatorname{sgn}}
\newcommand{\abs}[1]{\lvert#1\rvert}
\newcommand{\norm}[1]{\lVert#1\rVert}
\newcommand{\ltwonorm}[1]{\norm{#1}_{L^2(\Omega)}}
\newcommand{\ltwonormk}[1]{\norm{#1}_{L^2(K)}}
\newcommand{\linfnorm}[1]{\norm{#1}_{L^\infty(\Omega)}}
\newcommand{\linfnormk}[1]{\norm{#1}_{L^\infty(K)}}
\newcommand{\lh}{\abs{\ln h}}
\newcommand{\eps}{\varepsilon}
\newcommand{\T}{\mathcal{T}}
\newcommand{\lonenorm}[1]{\norm{#1}_{L^1(\Omega)}}
\renewcommand{\phi}{\varphi}
\newcommand{\lpnorm}[1]{\norm{#1}_{L^p(\Omega)}}
\newcommand{\Om}{\Omega}
\newcommand{\Abs}[1]{\left\lvert#1\right\rvert}
\newcommand{\mnorm}[1]{\norm{#1}_{\M(\Omega)}}

\section{Introduction}
In this paper we consider a semilinear elliptic equation of the following type
\begin{equation}\label{PDE:poisson_semilinear}
	\begin{aligned}
		-\Lap u(x) +d(x,u(x)) &= f(x) &\quad&\text{in } \Omega\\
		u(x) &=0 &\quad&\text{on } \partial \Omega
	\end{aligned}
\end{equation}
with a nonlinearity $d\colon \Omega \times \R \to \R$, posted on a convex polygonal/polyhedral domain $\Omega \subset \R^N$, $N=2,3$. Essentially we only assume that the nonlinearity $d$ is monotonically non-decreasing and continuous with respect to its second argument $u$, see below for precise assumptions.  This means, that non-differential, non-Lipschitz and even non-Hölder continuous nonlinearities are allowed. A typical example is a combustion type nonlinearity
\[
d(x,u) = \phi(x)r(u-\psi(x)) \quad \text{with }\; r(u) = \sgn(u) \abs{u}^s, \quad 0<s<1,
\]
and some $\phi,\psi \in L^\infty(\Omega)$ and $\phi(x) \ge 0$ for almost all $x \in \Omega$, see \cref{remark_assumptions} below. Such nonlinearities play a role by modeling of combustion problems and of flow problems in porous media, see, \cite{Nochetto:1988,KnabnerRannacher:2017} and the references therein. We refer also to \cite{Christof:2024} for recent results on optimal control of semilinear elliptic equations with such nonlinearities. 

For the finite element approximation $u_h$ in the space of linear finite elements $V_h$, see below for details, we will derive a priori error estimates for the error $u-u_h$. For the case of a Lipschitz continuous nonlinearity  error estimates are well known and have the same quality as for discretization of linear equations, i.e.,
\[
\ltwonorm{u-u_h} + h \ltwonorm{\nabla(u-u_h)} \le c h^2 \quad \text{and} \quad \linfnorm{u-u_h} \le ch^{2-\frac{N}{2}},
\]
where the constant $c$ depends on $\ltwonorm{f}$ and the nonlinearity $d$, see, e.g., \cite{CasasMateos:2002}. The case of a non-Lipschitz continuous nonlinearity is treated in \cite{Nochetto:1988} using a regularization approach. There, the domain $\Omega$ is assumed to be smooth and an $L^1$-$L^\infty$ duality argument is used to prove error estimates for $u_\eps- u_{\eps h}$ with respect to the $L^\infty(\Omega)$ norm, where $u_\eps$ is the solution of the equation with the regularized non-linearity and $u_{\eps h}$ is its finite element approximation. Similar approach is also applied in \cite{Barrett:1992} for a special Hölder continuous nonlinearity, see also \cite{BarrettShanahan:1991}. In \cite{KnabnerRannacher:2017} a parabolic problem with a Hölder continuous nonlinearity on a smooth domain is considered. There, also a variant of an $L^1$-$L^\infty$ duality argument is applied without any regularization to prove error estimates for a semi-discretization in space. We mention also the paper \cite{HafemeyerKahlePfefferer:2020}, where error estimates for the obstacle problem on a polygonal/polyhedral domain are shown. The authors prove there an interior pointwise estimate and then an $L^2$ estimate on strip along the boundary. This approach results in an optimal (up to some logarithmic terms) estimate with respect to the $L^2(\Omega)$ norm. We take a similar approach for the problem under consideration.

In the following we describe our contributions. For all results we do not require smoothness of the boundary of the domain $\Omega$. We assume throughout that $\Omega$ is polygonal for $N=2$, polyhedral for $N=3$ and convex. Moreover, we do not require any regularization of the nonlinearity.
\begin{itemize}
	\item Under general assumptions (see \cref{ass:1} below) on the nonlinearity $d$ (essentially monotonicity and continuity with respect to the second argument) we prove uniform boundedness of $\linfnorm{u_h}$, see \cref{FEM:semilinear:cor:uniform_bound_estimate} and the estimates
	\[
	\ltwonorm{\nabla(u-u_h)} \le ch \quad \text{and} \quad \linfnorm{u-u_h} \le ch^{2-\frac{N}{2}}
	\]
	for $f \in L^2(\Omega)$, see \cref{cor:linf_fl2} and \cref{theorem:h1_est}. For $f\in L^\infty(\Omega)$ the pointwise estimate is improved to
	\[
	\linfnorm{u-u_h} \le ch^{1+\gamma} \lh
	\]
	for every $0<\gamma<\gamma_\Omega$ with some $\gamma_\Omega \in (0,1]$ depending only on $\Omega$, see \cref{cor:linf_gamma}.
	\item Under an additional assumption (see \cref{ass:2} below), that the nonlinearity is Lipschitz-continuous in a neighborhood of boundary values, i.e., in our case in the interval $[-\rho,\rho]$ for some $\rho>0$, we prove
	\[
	\norm{u-u_h}_{L^\infty(\Omega_0)} \le ch^2 \lh^2 
	\]
	for an interior subdomain $\Omega_0$ with $\bar \Omega_0 \subset \Omega$ and $f \in L^\infty(\Omega)$, see \cref{theorem:interior}.
	 Moreover, we obtain the global estimate
	\[
	\ltwonorm{u-u_h} \le ch^2 \lh^2
	\]
	under the same assumptions, see \cref{theorem:L2}, which is our main result.
\end{itemize}

Throughout the paper we use standard notation for Lebesgue and Sobolev spaces.

The paper is structured as follows: In the following section we state precise assumptions on the nonlinearity $d$ and discuss properties of the solution of \cref{PDE:poisson_semilinear} as well as of an appropriate dual equation on the continuous level. In \Cref{sec:FEM} we recall finite element error estimates for the Ritz project under different levels of regularity.  In \Cref{sec:gen_er_est} we derive error estimates under the general assumption, i.e. under \cref{ass:1}. \Cref{sec:er_est_ass2} is devoted to the error estimates under the addition assumptions (see \cref{ass:2}), that the nonlinearity $d$ is Lipschitz continuous in a neighborhood of boundary values of the solution. In the last section we present a numerical example illustrating our results.

\section{Semilinear and dual equations}

In this section we specify assumptions on the nonlinearity $d$ and discuss existence, uniqueness and regularity for the state equation under consideration. Moreover, we discuss a dual equation to be used in the error analysis below.

\begin{assumption}\label{ass:1}
For the nonlinearity $d \colon \Omega \times \R \to \R$ we assume
\begin{itemize}
	\item For every $u_0 \in \R$, there holds $d(\cdot,u_0) \in L^\infty(\Omega)$
	\item For almost all $x\in \Omega$, the function $d(x,\cdot) \colon \R \to \R$ is monotonically non-decreasing
	\item For every $u_0 \in \R$ the nonlinearity $d$ is  equicontinuous at $u_0 \in \R$ in the following sense: For every $\eps >0$ there exists $\delta >0$ such that
\[
\abs{d(x,u)-d(x,u_0)} < \eps \quad \text{for all } u \in \R  \text{ with } \abs{u-u_0}< \delta  \text{ and } \text{for almost all } x \in \Omega.
\]
\end{itemize}
\end{assumption}

\begin{remark}\label{remark_assumptions}
For the class of nonlinearities given
by
\[
d(x,u) = \phi(x) r(u-\psi(x))
\]
with continuous and monotonically non-decreasing $r \colon \R \to \R$, $\phi,\psi \in L^\infty(\Omega)$ and $\phi(x) \ge 0$ for almost all $x \in \Omega$, \cref{ass:1} is fulfilled. To see this, one only has to discuss the equicontinuity. For given $u_0 \in \R$ consider the interval $I =
[u_0-2\linfnorm{\psi},u_0+2\linfnorm{\psi}]$. The function $r$ is continuous and therefore
uniformly continuous on $I$. Thus, for every $\eps >0$ there is $\delta >0$, which can be chosen
with $\delta < \linfnorm{\psi}$, such that
\[
\abs{r(\xi)-r(\eta)} < \frac{\eps}{\linfnorm{\phi}} \quad \text{for all } \xi,\eta \in I  \text{ with } \abs{\xi-\eta}<\delta.
\]
Let now $u \in \R$ with $\abs{u-u_0}<\delta$. Then
\[
u +\psi(x), u_0+ \psi(x) \in I \quad \text{and} \quad \abs{(u+\psi(x))-(u_0+\psi(x))} < \delta
\quad \text{for almost all } x \in \Omega.
\]
Therefore we obtain
\[
\abs{d(x,u)-d(x,u_0)} < \eps \quad \text{for almost all } x \in \Omega.
\]
\end{remark}

In the following, we will use the fact that a nonlinearity $d$ fulfilling
\cref{ass:1} is \emph{uniformly equicontinuous} on every bounded interval.
\begin{lemma}\label{FEM:lemma:uniform_equicontinuous}
Let $d \colon \Omega\times\R \to \R$ fulfill \cref{ass:1}. Then $d$, is uniformly equicontinuous on every compact interval $I \subset R$ in the following sense: For every $\eps >0$ there is $\delta >0$ such that 
\[
\abs{d(x,u)-d(x,v)} < \eps \quad \text{for all } u,v \in I \text{ with } \abs{u-v}< \delta
\text{ and almost all } x \in \Omega.
\]
\end{lemma}
\begin{proof}
The proof follows the lines of the classical proof that every continuous function is uniformly continuous on compact intervals.	
\end{proof}

For some of the results below we will require an additional assumption on the nonlinearity $d$.
\begin{assumption}\label{ass:2}
There are $\rho >0$ and $L>0$ such that for every $u,v \in [-\rho,\rho]$ there holds
\[
\abs{d(x,u)-d(x,v)} \le L \abs{u-v}
\]
for allmost all $x \in \Omega$.
\end{assumption}

\begin{remark}
The above assumption is fulfilled, e.g., for the nonlinearity
\[
d(x,u) = \phi(x)r(u-\psi(x)) \quad \text{with }\; r(u) = \sgn(u) \abs{u}^s, \quad 0<s<1,
\]
and some $\phi,\psi \in L^\infty(\Omega)$ and
\[ 
\phi(x) \ge 0, \; \psi(x) \ge \delta >0 
\]
for almost all $x \in \Omega$.
\end{remark}

\begin{remark}
\cref{ass:2} requires Lipschitz continuity in a neighborhood of boundary values. It can be correspondingly modified in the case of non-homogeneous Dirirchlet boundary conditions.
\end{remark}

In the sequel we will use also a modified (cut) nonlinearity $d_M \colon \Omega \times \R \to R$ with a cutting parameter $M>0$:
\begin{equation}\label{PDE:eq:cut_d}
	d_M(x,u) =  
	\begin{cases}
		d(x,-M), & \text{for } u < -M,\\
		d(x,u), & \text{for } {-M} \le u \le M,\\
		d(x,M), & \text{for } u > M.\\
	\end{cases}
\end{equation}

\begin{lemma}
Let $f \in L^p(\Omega)$ with some $p> \frac{N}{2}$ and let \cref{ass:1} be fulfilled. Then there is a unique weak solution $u_M \in H^1_0(\Omega) \cap C(\bar \Omega)$ fulfilling
\begin{equation}\label{eq_cut}
(\nabla u_M,\nabla \varphi) + (d_M(u_M),\varphi) = (f,\varphi) \quad \text{for all } \varphi \in H^1_0(\Omega).
\end{equation}
Moreover, there is a constant $c$ independent of $M$ such that
\[
\ltwonorm{\nabla u_M} + \linfnorm{u_M} \le c \left(\norm{f}_{L^p(\Omega)}+ \linfnorm{d(0)}\right).
\]
\end{lemma}
\begin{proof}
The uniqueness follows directly from monotonicity of $d$. The proof of existence of $u_M \in H^1_0(\Omega)$ is based on theory of monotone operators, see \cite[Theorem 26.A]{Zeidler:IIB:1990} and follows the lines of the proof of \cite[Theorem 4.4]{Troeltzsch:2010} (by changing Neumann to Dirichlet boundary conditions). The uniform (in $M$) boundedness of  $\linfnorm{u_M}$ is based on the Stampacchia methods \cite{KinderlehrerStampacchia:2000} and the ideas from \cite{Casas:1993}, see \cite[Theorem 4.7]{Troeltzsch:2010} for details.
\end{proof}

The uniform  (in $M$) boundedness of  $\linfnorm{u_M}$ leads directly to the existence of a solution with the original nonlinearity.

\begin{theorem}\label{theorem:u}
Let $f \in L^p(\Omega)$ with some $p> \frac{N}{2}$ and let \cref{ass:1} be fulfilled. Then there is a unique weak solution $u \in H^1_0(\Omega) \cap C(\bar \Omega)$ fulfilling
\[
(\nabla u,\nabla \varphi) + (d(u),\varphi) = (f,\varphi) \quad \text{for all } \varphi \in H^1_0(\Omega).
\]
Moreover, there is a constant $c$ such that
\[
\ltwonorm{\nabla u} + \linfnorm{u} \le c \left(\norm{f}_{L^p(\Omega)}+ \linfnorm{d(0)}\right).
\]
If $f \in L^2(\Omega)$ then $u \in H^2(\Omega)$ and the following estimate holds
\[
\ltwonorm{\nabla^2 u} \le c \left(\ltwonorm{f}+ \linfnorm{d(M)}+ \linfnorm{d(-M)}\right),
\]
where
\[
M = c \left(\norm{f}_{L^p(\Omega)} + \linfnorm{d(0)}\right)
\]
and the constant $c$ is independent of $f$ and $d$.
\end{theorem}
\begin{proof}
The existence, uniqueness and the first estimate follow directly from the previous lemma. Then a bootstrapping argument and $H^2$ regularity on convex domains for the linear Poisson equation lead to the corresponding estimate for $\ltwonorm{\nabla^2 u}$.
\end{proof}

For the numerical analysis below we require the dual equation
\begin{equation}\label{eq:dual}
	\begin{aligned}
		-\Lap z(x) +b(x)z(x) &= \psi(x) &\quad&\text{in } \Omega\\
		z(x) &=0 &\quad&\text{on } \partial \Omega
	\end{aligned}
\end{equation}
for given $b \ge 0$ and $\psi$ to be specified later. The following
result, cf. \cite[Lemma 2.1]{Nochetto:1988}, provides an estimate with a constant independent of $b$.

\begin{lemma}\label{FEM:semilinear:lemma:L1_est_b}
	Let $\psi \in L^2(\Omega)$ and $b \in L^\infty(\Omega)$ with $b \ge 0$ almost everywhere on
	$\Omega$. Then, the solution $z \in H^1_0(\Omega) \cap H^2(\Omega)$ of
	\begin{equation}\label{FEM:eq:in_lemma:L1_est_b}
		(\nabla z,\nabla \phi) + (bz,\phi) = (\psi,\phi) \quad \text{for all } \phi \in H^1_0(\Omega)
	\end{equation}
	fulfills
	\[
	\ltwonorm{z} + \lonenorm{\Lap z} + \lonenorm{bz} \le c\lonenorm{\psi},
	\]
	where constant $c$ on the right-hand side of the estimate is independent of $b$ and $\psi$.
\end{lemma}

\section{Finite element discretization}\label{sec:FEM}
Let $\set{\T_h}$ be a family of shape-regular, quasi-uniform meshes. Let $V_h \subset H^1_0(\Omega)\cap C(\bar \Omega)$ be the space linear finite elements on $\T_h$. We will use Lagrange interpolation operator $i_h \colon C_0(\Omega) \to V_h$, which fulfills standard interpolation estimates. Moreover, we will use the Ritz projection $R_h \colon W^{1,1}_0(\Omega) \to V_h$ defined by
\begin{equation}\label{eq:Ritz}
(\nabla(u-R_h u), \nabla \phi_h)  = 0 \quad \text{for all } \phi_h \in V_h.
\end{equation}
We will use the following estimates for $u-R_h u$.
\begin{lemma}\label{lemma:standard_Rh}
Let $u \in H^2(\Omega)\cap H^1_0(\Omega)$. Then there holds
\[
\ltwonorm{u-R_h u} + h \ltwonorm{\nabla(u-R_h u)} \le c h^2 \ltwonorm{\Lap u}
\]
and
\[
\linfnorm{u-R_h u} \le ch^{2-\frac{N}{2}}\ltwonorm{\Lap u}.
\]
\end{lemma}
\begin{proof}
These results are classical, see, e.g., \cite{BrennerScott:2008}.
\end{proof}

The following improved pointwise estimate holds provided $\Lap u \in L^\infty(\Omega)$.
\begin{lemma}\label{lemma:gl_est_linf}
There is $\gamma_\Omega \in (0,1]$ depending only on the domain $\Omega$ and a constant $c>0$ such that for every $0<\gamma<\gamma_\Omega$ and every $u \in H^1_0(\Omega)$ with $\Lap u \in L^\infty(\Omega)$ there holds
\[
\linfnorm{u-R_h u} \le ch^{1+\gamma}\lh\linfnorm{\Lap u}.
\]
\end{lemma}
\begin{proof}
First we discuss the Hölder continuity of the derivatives of the solution $u$ to a linear elliptic equation
\[
	\begin{aligned}
		-\Lap u(x) &= f(x) &\quad&\text{in } \Omega\\
		u(x) &=0 &\quad&\text{on } \partial \Omega
	\end{aligned}
\]
with $f = -\Lap u \in L^\infty(\Omega)$. For $N=2$ we can use convexity of $\Omega$ and argue by $W^{2,p}$ regularity. By \cite[Theorem 4.4.3.7]{Grisvard:1985} there exists $p_\Omega > 2$ such that $u \in W^{2,p}(\Omega)$ holds for every $2 \le p < p_\Omega$. By a classical embedding theorem we that obtain the existence of $0<\gamma_\Omega\le 1$ such that $u \in C^{1,\gamma}(\bar \Omega)$ holds for every $0< \gamma<\gamma_\Omega$ with
\[
\norm{u}_{C^{1,\gamma}(\bar \Omega)} \le c \linfnorm{f}.
\]

For $N=3$ the corresponding result follows from Hölder estimates of derivatives of the Green's function on a convex polyhedral domain, see \cite[Theorem 1]{GuzmanLeykekhmanannRossmannSchatz:2009}.

As the next step we use a best approximation (up to a logarithmic term) result  on convex domains for the error in $L^\infty(\Omega)$, i.e., 
\[
\linfnorm{u-R_h u} \le c \lh \inf_{v_h \in V_h}\linfnorm{u-v_h},
\]
see \cite{FrehseRannacher:1976} for $N=2$ and \cite{LeykekhmanD_VexlerB_2016c} for $N=3$. We choose $v_h = i_h u$ and obtain by standard interpolation estimates
\[
\linfnorm{u-R_h u} \le c \lh \linfnorm{u-i_h u} \le c h^{1+\gamma} \lh \norm{u}_{C^{1,\gamma}(\bar \Omega)} \le  c h^{1+\gamma} \lh \linfnorm{f}.
\]
This completes the proof.
\end{proof}

\begin{remark}
For $N=2$, we have 
\[
\gamma_\Omega = 
\begin{cases}
	\lambda_\Omega - 1 & \text{for } \lambda_\Omega < 2,\\
	1 & \text{for } \lambda_\Omega \ge 2,
\end{cases}
\]
where $\lambda_\Om = \nicefrac{\pi}{\omega_\Omega}$ and $\omega_\Omega$ is the largest interior angle of $\Omega$. For $N=3$, the value of $\gamma_\Omega$ depends on both the geometry of vertices and edges of
$\Omega$, see discussion in \cite{GuzmanLeykekhmanannRossmannSchatz:2009}. For $\Omega =(0,1)^3\subset
\R^3$ being a cube, there holds $\gamma_\Omega=1$. 
\end{remark}

The following result provides an interior estimate.
\begin{lemma}\label{lemma:interior_est_Rh}
Let $u \in H^1_0(\Omega)$. Let $\Omega_0$ and $\Omega_1$ be interior subdomains with $\bar \Omega_0 \subset \Omega_1$ and $\bar \Omega_1 \subset \Omega$. Let $\Lap u \in L^2(\Omega)$ and $\Lap u \in L^\infty(\Omega_1)$. Then there holds
\[
\norm{u-R_h u}_{L^\infty(\Omega_0)} \le c h^2 \lh^2 \left(\norm{\Lap u}_{L^\infty(\Omega_1)}+\ltwonorm{\Lap u}\right).
\]
\end{lemma}
\begin{proof}
We consider an intermediate subdomain $\Omega_0'$ with $\bar \Omega_0 \subset \Omega_0'$ and $\bar \Omega_0' \subset \Omega_1$. By the interior regularity result we obtain $u \in W^{2,p}(\Omega_0')$ for every $2\le p <\infty$ with
\[
\norm{u}_{W^{2,p}(\Omega_0')} \le c p \left(\norm{f}_{L^\infty(\Omega_1)} +  \ltwonorm{f}\right),
\]
where $f= -\Lap u$ and a constant $c$ is independent of $u$ and $p$. The explicit dependence on $p$ in the previous estimate can be traced from the proof of \cite[Theorem 9.9]{GilbargTrudinger:2001}.

By the interior estimate from \cite[Theorem 5.1]{SchatzWahlbin:1977} we get for every such $p$
\[
\begin{aligned}
	\norm{u-u_h}_{L^\infty(\Omega_0)} &\le c \lh \norm{u-i_h u}_{L^\infty(\Omega_0')} + \ltwonorm{u-u_h}\\
	&\le c h^{2-\frac{N}{p}} \lh \norm{\nabla^2 u}_{L^p(\Omega_0')}+ \ltwonorm{u-u_h}\\
	&\le c p  h^{2-\frac{N}{p}} \lh \left(\norm{f}_{L^\infty(\Omega_1)} +  \ltwonorm{f}\right) + \ltwonorm{u-u_h},
\end{aligned}
\]
where we used standard interpolation estimates. We choose $p=\lh$ and obtain using
\[
h^{-\frac{N}{p}} = h^{-\frac{N}{\lh}} = h^{\frac{N}{\ln h}} = e^N
\]
the estimate
\[
\norm{u-u_h}_{L^\infty(\Omega_0)} \le c h^2 \lh^2 \left(\norm{f}_{L^\infty(\Omega_1)} +  \ltwonorm{f}\right) + \ltwonorm{u-u_h},
\]
which completes the proof by the estimate for $\ltwonorm{u-u_h}$ from \cref{lemma:standard_Rh}.
\end{proof}

We require also error estimates for the Ritz projection of the solution with low regularity, especially if $\Lap u \in L^1(\Omega)$.

\begin{lemma}\label{lemma:est_l2_Lap_l1}
Let $u \in W^{1,1}_0(\Omega)$ and $\Lap u \in L^1(\Omega)$. Then there holds
\[
\ltwonorm{u-R_h u} \le ch^{2-\frac{N}{2}}\lonenorm{\Lap u}.
\]
\end{lemma}
\begin{proof}
The straightforward proof uses a duality argument and the pointwise estimate from \cref{lemma:standard_Rh}.
\end{proof}

\begin{lemma}\label{lemma:est_l1_Lap_l1}
	Let $\gamma_\Omega \in (0,1]$ be as in \cref{lemma:gl_est_linf} and let $0<\gamma<\gamma_\Omega$. Let $u \in W^{1,1}_0(\Omega)$ and $\Lap u \in L^1(\Omega)$. Then there holds
	\[
	\lonenorm{u-R_h u} \le ch^{1+\gamma}\lh\lonenorm{\Lap u}.
	\]
\end{lemma}
\begin{proof}
	The straightforward proof uses a duality argument and the pointwise estimate from \cref{lemma:gl_est_linf}.
\end{proof}

\section{General error estimates}\label{sec:gen_er_est}

The discrete approximation to \cref{PDE:poisson_semilinear} is defined by
\begin{equation}\label{FEM:poisson_semilinear_weak_h}
u_h \in V_h \quad:\quad (\nabla u_h, \nabla \phi_h) + (d(u_h),\phi_h)= (f,\phi_h) \quad \text{for all } \phi_h \in V_h.
\end{equation}

\begin{theorem}\label{FEM:SemiLinear_Existence_h}
	Let $f \in L^p(\Omega)$ with some $p> \frac{N}{2}$ and let \cref{ass:1} be fulfilled. Then, there exists a
	unique solution of \cref{FEM:poisson_semilinear_weak_h}. There is a constant $c$ independent of
	$f$, $d$, and $h$ such that
	\[
	\ltwonorm{\nabla u_h} \le c \left(\lpnorm{f} + \linfnorm{d(0)}\right).
	\]
\end{theorem}
\begin{proof}
The proof is based on Brouwer's fixed point theorem, see \cite[Chapter 1, Satz 2.23]{Ruzicka:2020} for the idea of the proof.
\end{proof}

For given $M>0$ we define an approximation $u_{M,h} \in V_h$ of the solution to the equation with cut nonlinearity \cref{eq_cut} the same way as \cref{FEM:poisson_semilinear_weak_h} replacing $d$ by $d_M$. We obtain the following Galerkin orthogonality relations:
\begin{equation}\label{FEM:semilinear_Galerkin_orthogonality}
	(\nabla (u-u_{h}),\nabla \phi_h) + ( d(u)-d(u_{h}),\phi_h ) = 0 \quad \text{for all } \phi_h \in V_h
\end{equation}
and
\begin{equation}\label{FEM:semilinear_Galerkin_orthogonality_M}
	(\nabla (u_M-u_{M,h}),\nabla \phi_h) + ( d_M(u_M)-d_M(u_{M,h}),\phi_h ) = 0 \quad \text{for all } \phi_h \in V_h.
\end{equation}

The next theorem provides a pointwise estimate for $u_M - u_{M,h}$.

\begin{theorem}\label{FEM:semilinear:theorem_linf_est_u_M}
	Let the nonlinearity $d$ fulfill \cref{ass:1}. Let
	$f \in L^2(\Omega)$.
	Let $u_M \in
	H^1_0(\Omega)\cap C(\bar \Omega)$ be the solution of \cref{eq_cut} and
	$u_{M,h} \in V_h$ its discrete approximation. Then there holds
	\[
	\linfnorm{u_M-u_{M,h}} \le c_M h^{2-\frac{N}{2}} \left(\ltwonorm{f}+1\right),
	\]
	where the constant $c_M$ may depend on $M$ but is independent of $h$, $f$, and $d$.
\end{theorem}
\begin{proof}
	We denote $e=u_M - u_{M,h}$.
	
	\emph{Step 1: Construction of an adjoint problem.} For a duality argument, we will consider the adjoint problem
	\begin{equation}\label{FEM:semilinear:eq:adjoint_in_linf_proof}
		z \in H^1_0(\Omega) \quad \colon \quad (\nabla z,\nabla v) + (bz,v) = (\psi,v) \quad \text{for all } v \in H^1_0(\Omega)
	\end{equation}
	for some $\psi \in L^2(\Omega)$ with $\lonenorm{\psi}\le 1$. Since the nonlinearity $d$ is, in general, not differentiable, we can not use the derivative of $d$ for the definition of $b$. Also, a
	differential quotient can not be used directly since $d$ is not necessarily Lipschitz continuous.
	The following definition of the coefficient $b$ is inspired by \cite[p. 188]{KnabnerRannacher:2017}. We first define 
	\[
	R_{M,h} = \linfnorm{u_{M,h}} + \linfnorm{u_M}.
	\]
	Since $d$ (and therefore also $d_M$) fulfills \cref{ass:1}, we obtain by \cref{FEM:lemma:uniform_equicontinuous} that $d_M$ is uniformly equicontinuous on the interval $I_{M,h} = [-R_{M,h},R_{M,h}]$. Therefore, there exists $\rho_{M,h}>0$ such that for all $\xi,\eta \in I_{M,h}$ with $\abs{\xi-\eta}< \rho_{M,h}$ there holds 
	\begin{equation}\label{FEM:semilinear:appl_equicont}
		\abs{d_M(x,\xi)-d_M(x,\eta)} < h^2 \quad \text{for almost all } x \in \Omega.
	\end{equation}
	Using this, we define the coefficient function $b$ for almost all $x\in\Om$ by
	\begin{equation}\label{FEM:semilinear:def_b}
		b(x) =
		\frac{d_M(x,u_M(x))-d_M(x,u_{M,h}(x))}{\sgn(e(x))\max(\abs{e(x)},\rho_{M,h})}.
	\end{equation}
	By this definition, we directly obtain $b \in L^\infty(\Omega)$, but $\linfnorm{b}$ is, in general, not bounded uniformly in $h$. Moreover, we have by the monotonicity of $d_M$ that
	\[
	b(x) \ge 0 \quad \text{for almost all } x \in \Omega.
	\]
	In addition, we prove the following property of $b$: There holds
	\[
	b(x)e(x) = \left(d(x,u_M(x))-d(x,u_{M,h}(x))\right) + e_b(x) \quad \text{for almost all } x \in \Omega
	\]
	with $e_b \in L^\infty(\Omega)$ and $\linfnorm{e_b} \le 2 h^2$.
	To see this, we distinguish two cases:
	In the case $\abs{e(x)} \ge \rho_{M,h}$, the denominator in \cref{FEM:semilinear:def_b} is equal
	to $e(x)$. Therefore, it holds
	\[
	b(x)e(x) = d(x,u_M(x))-d(x,u_{M,h}(x))
	\]
	and consequently $e_b(x)=0$.
	If $\abs{e(x)} < \rho_{M,h}$, we have by \cref{FEM:semilinear:appl_equicont}
	\[
	\begin{aligned}
		\abs{e_b(x)} &= \Abs{b(x)e(x) - \left(d(x,u_M(x))-d(x,u_{M,h}(x))\right)}\\
		&=\Abs{d(x,u_M(x))-d(x,u_{M,h}(x))} \Abs{\frac{e(x)}{\sgn(e(x)) \rho_{M,h}} - 1}\\
		&\le h^2 \left(\frac{\abs{e(x)}}{\rho_{M,h}}+1\right) \le 2 h^2.
	\end{aligned}
	\]
	
	\emph{Step 2: Error estimation.}
	To estimate $\linfnorm{e}$, we first observe that
	\[
	\linfnorm{e} = \sup_{\substack{\psi \in L^2(\Omega)\\\lonenorm{\psi}\le 1}} (e,\psi).
	\]
	Let $\psi \in L^2(\Omega)$ with
	$\lonenorm{\psi}\le 1$. We consider the corresponding solution $z \in H^1_0(\Omega)$ of
	\cref{FEM:semilinear:eq:adjoint_in_linf_proof} with the coefficient $b$ defined by
	\cref{FEM:semilinear:def_b}. By choosing  $\phi=e$ in
	\cref{FEM:semilinear:eq:adjoint_in_linf_proof}, using properties of $b$ as well as
	\cref{FEM:semilinear_Galerkin_orthogonality_M} with $\phi_h$ chosen as Ritz projection $\phi_h = R_h
	z$, we get
	\[
	\begin{aligned}
		(e,\psi) 
		&= (\nabla e,\nabla z) + (be,z)\\
		&= (\nabla e,\nabla z) + ( d_M(u_M)-d_M(u_{M,h}),z ) + (e_b,z)\\
		&= (\nabla e,\nabla (z-R_h z)) + ( d_M(u_M)-d_M(u_{M,h}),z-R_h z) + (e_b,z)\\
		&= (\nabla u_M,\nabla (z-R_h z)) + ( d_M(u_M)-d_M(u_{M,h}),z-R_h z) + (e_b,z),
	\end{aligned}
	\]
	where in the last step, the definition \cref{eq:Ritz} of the Ritz projection was
	employed. Using now equation \cref{eq_cut}, we obtain
	\[
	\begin{aligned}
		(e,\psi) 
		&= (\nabla u_M,\nabla (z-R_h z)) + ( d_M(u_M)-d_M(u_{M,h}),z-R_h z) + (e_b,z)\\
		&= (f-d_M(u_M),z-R_h z)+ ( d_M(u_M)-d_M(u_{M,h}),z-R_h z) + (e_b,z)\\
		&= (f-d_M(u_{M,h}),z-R_h z) + (e_b,z)\\
		& \le \left(\ltwonorm{f} + c\linfnorm{d_M(u_{M,h})}\right) \ltwonorm{z-R_h z}\\
		&\quad+ \linfnorm{e_b} \lonenorm{z}.
	\end{aligned}
	\] 
	We apply 
	\cref{lemma:est_l2_Lap_l1}
	leading to
	\begin{equation}\label{FEM:eq:in:proof:r-Rhz_lpp}
		\begin{aligned}
			\ltwonorm{z-R_h z} &\le  c h^{2-\frac{N}{2}}\lonenorm{\Lap z }\\
			& \le c h^{2-\frac{N}{2}} \lonenorm{\psi}  = ch^{2-\frac{N}{2}},
		\end{aligned}
	\end{equation}
	where we have used the fact that $\Lap z \in L^1(\Omega)$ and the estimate from
	\cref{FEM:semilinear:lemma:L1_est_b}. For $\lonenorm{z}$, we estimate also using \cref{FEM:semilinear:lemma:L1_est_b} 
	\[
	\lonenorm{z} \le c \ltwonorm{z}  \le c \lonenorm{\psi}  = c.
	\]
	Using the definition \cref{PDE:eq:cut_d} of $d_M$, we have
	\[
	\linfnorm{d_M(u_{M,h})} \le \linfnorm{d(M)} + \linfnorm{d(-M)} \le c_M.
	\]
	Putting all estimates together, we obtain
	\[
	(e,\psi) \le ch^{2-\frac{N}{2}} \left(\ltwonorm{f} + c_M\right) + ch^2.
	\]
	This completes the proof.
\end{proof}

We obtain the uniform boundedness of $\linfnorm{u_h}$  as a direct corollary.

\begin{corollary}\label{FEM:semilinear:cor:uniform_bound_estimate}
	Let the nonlinearity $d$ fulfill \cref{ass:1}. Let 
$f \in L^2(\Omega)$  and $u \in H^1_0(\Omega) \cap C(\bar \Omega)$ be the solution of \cref{PDE:poisson_semilinear}. Let $u_h \in V_h$ be the solution of \cref{FEM:poisson_semilinear_weak_h}. There exists $h_0>0$ such that for all $h<h_0$ there holds
\[
\linfnorm{u_h} \le 2\linfnorm{u}.
\]
\end{corollary}
\begin{proof}
		Without loss of generality, we can assume $u \neq 0$. Indeed, if $u=0$ then we have $u_h=0$
	by \cref{FEM:semilinear_Galerkin_orthogonality} and the desired inequality holds. 
	We choose $M = 2\linfnorm{u} >0$. Then, there holds for the unique solution $u_M$ of
	\cref{eq_cut} that $u_M=u$ by construction of the cut nonlinearity $d_M$
	\cref{PDE:eq:cut_d}. For the discrete solution $u_{M,h} \in V_h$ we have by \cref{FEM:semilinear:theorem_linf_est_u_M} that
	\[
	\linfnorm{u-u_{M,h}} = \linfnorm{u_M-u_{M,h}} \le c_M h^{2-\frac{N}{2}} \left(\lpnorm{f}+1\right).
	\]
	We choose $h_0>0$ such that 
	\[
	c_M h_0^{2-\frac{N}{2}} \left(\lpnorm{f}+1\right) < \linfnorm{u}.
	\]
	Then, we have for all $0<h<h_0$
	\[
	\linfnorm{u_{M,h}} \le \linfnorm{u} +  c_M h^{2-\frac{N}{2}} \left(\lpnorm{f}+1\right) <2 \linfnorm{u} = M.
	\]
	Thus, there we have $u_{M,h} = u_h$ and the estimate
	\[
	\linfnorm{u_h} \le 2\linfnorm{u}
	\]
	holds.
\end{proof}

\begin{remark}
	The statement of the above corollary holds under a weaker assumption on $f$, i.e. $f \in L^p(\Omega)$ with some $p> \frac{N}{2}$. Thus, one obtains the same estimate for $\linfnorm{u_h}$ as on the continuous level, i.e.
	\[
	\linfnorm{u_h} \le c \left(\norm{f}_{L^p(\Omega)}+ \linfnorm{d(0)}\right),
	\]
	cf. \cref{theorem:u}.
\end{remark}

The following two corollaries provide $L^\infty(\Omega)$ error estimates for $u-u_h$ under different assumptions an $f$.
\begin{corollary}\label{cor:linf_fl2}
	Let the nonlinearity $d$ fulfill the \cref{ass:1}.
	Let $f \in L^2(\Omega)$ and $u \in H^1_0(\Omega)\cap C_0(\Omega)$ be the solution of
	\cref{PDE:poisson_semilinear}. Let $u_h \in V_h$ be  the solution of
	\cref{FEM:poisson_semilinear_weak_h}. Then, there holds
	\[
	\linfnorm{u-u_h} \le c h^{2-\frac{N}{2}}.
	\]
	The constant $c$ depends on the nonlinearity $d$ and $\ltwonorm{f}$ and is independent of $h$.
\end{corollary}
\begin{proof}
	The estimate follows directly from \cref{FEM:semilinear:theorem_linf_est_u_M} since $u_M=u$ and $u_{M,h}=u_h$ by the appropriate choice of $M$ as in \cref{FEM:semilinear:cor:uniform_bound_estimate}. 
\end{proof}

\begin{corollary}\label{cor:linf_gamma}
		Let $\gamma_\Omega \in (0,1]$ be as in \cref{lemma:gl_est_linf} and let $0<\gamma<\gamma_\Omega$. Let the nonlinearity $d$ fulfill the \cref{ass:1}.
	Let $f \in L^\infty(\Omega)$ and $u \in H^1_0(\Omega)\cap C_0(\Omega)$ be the solution of
	\cref{PDE:poisson_semilinear}. Let $u_h \in V_h$ be  the solution of
	\cref{FEM:poisson_semilinear_weak_h}. Then, there holds
	\[
	\linfnorm{u-u_h} \le c h^{1+\gamma}\lh.
	\]
	The constant $c$ depends on the nonlinearity $d$ and $\linfnorm{f}$ and is independent of $h$.
\end{corollary}
\begin{proof}
Inspecting the proof of \cref{FEM:semilinear:theorem_linf_est_u_M} and setting $M = 2 \linfnorm{u}$, such that $u_M = u$ and $u_{M,h} = u_h$, we obtain for $e = u-u_h$ and $\psi \in L^2(\Omega)$ with $\lonenorm{\psi} \le 1$
\[
\begin{aligned}
(e,\psi) &= (f-d(u_h),z-R_h z) + (e_b,z)\\
&\le \left(\linfnorm{f} + c\linfnorm{d(u_h)}\right) \lonenorm{z-R_h z}+ \linfnorm{e_b} \lonenorm{z}.
\end{aligned}
\]
We apply \cref{lemma:est_l1_Lap_l1} leading to
\[
\lonenorm{z-R_h z} \le c h^{1+\gamma} \lh \lonenorm{\Lap z},
\]
the fact that $\linfnorm{e_b} \le ch^2$ by construction and obtain
\[
(e,\psi) \le c h^{1+\gamma} \lh\lonenorm{\Lap z} + ch^2 \lonenorm{z} \le c h^{1+\gamma} \lh,
\]
where we used \cref{FEM:semilinear:lemma:L1_est_b} in the last step. This completes the proof.
\end{proof}

\begin{remark}
In the case of the unit square $\Omega = (0,1)^2$ or unit cube $\Omega = (0,1)^3$ we have $\gamma_\Omega = 1$ and thus we get for $\linfnorm{u-u_h}$ almost ${\mathcal O}(h^2)$ estimate in these cases.
\end{remark}

The following theorem provides an optimal estimate in the energy norm.

\begin{theorem}\label{theorem:h1_est}
		Let the nonlinearity $d$ fulfill the \cref{ass:1}.
	Let $f \in L^2(\Omega)$ and $u \in H^1_0(\Omega)\cap C_0(\Omega)$ be the solution of
	\cref{PDE:poisson_semilinear}. Let $u_h \in V_h$ be  the solution of
	\cref{FEM:poisson_semilinear_weak_h}. Then, there holds
	\[
	\ltwonorm{\nabla(u-u_h)} \le c h.
	\]
	The constant $c$ depends on the nonlinearity $d$ and $\ltwonorm{f}$ and is independent of $h$.
\end{theorem}	
\begin{proof}
We denote the error by $e= u-u_h$ and decompose it as $e = \eta + \xi_h$ with
\[
\eta = u - R_h u \quad \text{and} \quad \xi_h = R_h u - u_h \in V_h,
\]
where $R_h u$ is the Ritz projection of $u$. There holds by the orthogonality relation \cref{FEM:semilinear_Galerkin_orthogonality}
\[
\begin{aligned}
	\ltwonorm{\nabla e}^2 &= (\nabla e,\nabla \eta) + (\nabla e,\nabla \xi_h)\\
	&= (\nabla e,\nabla \eta)-(d(u)-d(u_h),\xi_h)\\
	&= (\nabla e,\nabla \eta)-(d(u)-d(u_h),e) + (d(u)-d(u_h),\eta)\\
	&\le (\nabla e,\nabla \eta) + (d(u)-d(u_h),\eta),
\end{aligned}
\]
where in the last step, we have used the monotonicity of $d$. By the definition of $R_h$, we get
\[
(\nabla e,\nabla \eta) = (\nabla u, \nabla \eta) = (f-d(u),\eta)
\]
resulting in
\[
\ltwonorm{\nabla e}^2 \le (f-d(u_h),\eta).
\]
By the boundedness of $u_h$ from \cref{FEM:semilinear:cor:uniform_bound_estimate} we get
\[
\ltwonorm{\nabla e}^2 \le c \ltwonorm{\eta}.
\]
The estimate for the Ritz projection from \cref{lemma:standard_Rh} leads to
\[
\ltwonorm{\nabla e}^2 \le c h^2 \ltwonorm{\Lap u}.
\]
Using
\[
\ltwonorm{\Lap u} \le \ltwonorm{f} + \ltwonorm{d(u)} \le  c
\]
completes the proof.
\end{proof}

\section{Error estimates under \cref{ass:2}}\label{sec:er_est_ass2}
For given $\delta >0$ we define a boundary strip 
\[
D_\delta = \Set{x \in \Omega| \operatorname{dist}(x,\partial \Omega) < \delta}.
\]
Since the solution $u$ of \cref{PDE:poisson_semilinear} fulfills $u \in C_0(\Omega)$ there is $\delta >0$, such that
\begin{equation}\label{eq:D_4d}
\abs{u(x)} \le \frac{\rho}{2} \quad \text{for all }\; x \in D_{4\delta}
\end{equation}
with $\rho>0$ from \cref{ass:2}.

\begin{theorem}\label{theorem:interior}
Let the nonlinearity $d$ fulfill the \cref{ass:1} and \cref{ass:2}. Let $\Omega_0$ be an interior subdomain, i.e. $\bar \Omega_0 \subset \Omega$. Let $f \in L^\infty(\Omega)$ and $u \in H^1_0(\Omega)\cap C_0(\Omega)$ be the solution of
\cref{PDE:poisson_semilinear}. Let $u_h \in V_h$ be  the solution of
\cref{FEM:poisson_semilinear_weak_h}. Then, there holds
\[
\norm{u-u_h}_{L^\infty(\Omega_0)} \le c h^2 \lh^2.
\]
\end{theorem}
\begin{proof}
As before we denote $e = u-u_h$. There holds
\[
\norm{e}_{L^\infty(\Omega_0)} = \sup\Set{(e,\psi) | \psi \in L^2(\Omega), \lonenorm{\psi}\le 1, \operatorname{supp} \psi \subset \Omega_0}.
\]
For such a function $\psi$ we consider the dual equation \cref{FEM:semilinear:eq:adjoint_in_linf_proof} with $b$ defined as in the proof of \cref{FEM:semilinear:theorem_linf_est_u_M} by \cref{FEM:semilinear:def_b}. Please note, that $M$ is chosen as $M = 2 \linfnorm{u}$ such that $u_M = u$ and $u_{M,h} = u_h$, cf. the proof of \cref{cor:linf_gamma}. As there, we obtain
\[
(e,\psi) \le \left(\linfnorm{f} + c\linfnorm{d(u_h)}\right) \lonenorm{z-R_h z}+ \linfnorm{e_b} \lonenorm{z}.
\]
Recalling that $\linfnorm{e_b} \le ch^2$, it remains to estimate $\lonenorm{z-R_h z}$. To this end we consider 
\[
y \in H^1_0(\Omega) \quad:\quad (\nabla y,\nabla \phi) = (g,\phi) \quad \text{for all }\; \phi \in H^1_0(\Omega)
\]
with $g = \operatorname{sgn}(z-R_h z)$. Note, that 
\[
\linfnorm{\Lap y} = \linfnorm{g} \le 1.
\]
There holds
\[
\begin{aligned}
\lonenorm{z-R_h z} &= (z-R_h z, g) = (\nabla y,\nabla z) - (\nabla y,\nabla R_h z)\\
& = (\nabla y,\nabla z) - (\nabla R_h y,\nabla z) = (\nabla(y-R_h y),\nabla z)\\
& = (\psi,y-R_h y) - (bz,y-R_h y),
\end{aligned}
\]
where we have used the definition of the Ritz projection and the equation \cref{FEM:semilinear:eq:adjoint_in_linf_proof} for $z$. We get by $\operatorname{supp} \psi \subset \Omega_0$
\[
(\psi,y-R_h y) \le \lonenorm{\psi} \norm{y-R_h y}_{L^\infty(\Omega_0)} \le c h^2\lh^2 \linfnorm{g} \le c h^2\lh^2
\]
by \cref{lemma:interior_est_Rh}. To estimate $(bz,y-R_h y)$ we observe that by the definition \cref{FEM:semilinear:def_b} of $b$, by \cref{eq:D_4d} and by \cref{ass:2} we have
\[
0 \le b(x) \le \frac{d(x,u(x))-d(x,u_h(x))}{e(x)} \le L \quad \text{for almost all }\; x \in D_{4\delta}.
\]
Thus, we have
\[
\norm{b z}_{L^2(D_{4\delta})} \le L \ltwonorm{z} \le c.
\]
This allows us to estimate
\[
\begin{aligned}
(bz,y-R_h y) &\le \norm{b z}_{L^2(D_{4\delta})} \norm{y-R_h y}_{L^2(D_{4\delta})} + \norm{b z}_{L^1(\Omega \setminus D_{4\delta})} \norm{y-R_h y}_{L^\infty(\Omega \setminus D_{4\delta})}\\
& \le c \ltwonorm{y-R_h y} + c \norm{y-R_h y}_{L^\infty(\Omega \setminus D_{4\delta})}.
\end{aligned}
\]
Using
\[
\ltwonorm{y-R_h y} \le c h^2 \ltwonorm{\Lap y} = c h^2 \ltwonorm{g} \le ch^2
\]
from \cref{lemma:standard_Rh} and
\[
\norm{y-R_h y}_{L^\infty(\Omega \setminus D_{4\delta})} \le ch^2 \lh^2 \linfnorm{g} \le ch^2 \lh^2
\]
from interior estimate in \cref{lemma:interior_est_Rh} we obtain
\[
(bz,y-R_h y) \le ch^2 \lh^2.
\]
Putting terms together we complete the proof.
\end{proof}

The above result provides an error estimate on an interior domain with respect to the $L^\infty(\Omega_0)$ norm. Since $\Omega$ is a general convex (polygonal/polyhedral) domain, we can not expect an estimate of this quality ($\mathcal O(h^2)$ up to logarithmic terms) with respect to $L^\infty(\Omega)$ norm, i.e. on the entire domain $\Omega$. Thus, our goal is to prove an estimate of this order with respect to $L^2(\Omega)$ norm. To this end we estimate the error on the boundary strip $D_\delta$. We define a smooth cut-off function $\omega \colon \bar \Omega \to [0,1]$ with
\begin{equation}\label{eq:omega}
\omega(x)
\begin{cases}
=1, & x \in D_{2\delta}\\
\in [0,1], & x \in D_{3\delta} \setminus D_{2\delta}\\
=0 & x \in \Omega \setminus D_{3\delta}
\end{cases}
\end{equation}
with $\abs{\nabla \omega(x)} \le \delta^{-1}$ and  $\abs{\nabla^2 \omega(x)} \le \delta^{-2}$.

We will use the following superapproximation estimates. Such results can be found in several publications starting from \cite{NitscheSchatz:1974}. We refer to \cite{HafemeyerKahlePfefferer:2020} for a result very close to the following, where the terms on the right-hand side are supported on $\Omega \setminus D_\delta$.
\begin{lemma}\label{lemma:superconv}
Let $\omega$ by defined by \cref{eq:omega} and let $h< \delta$. Then there holds
\begin{multline*}
\ltwonorm{(\operatorname{id} - i_h)(\omega^2 v_h)} + h \ltwonorm{\nabla(\operatorname{id} - i_h)(\omega^2 v_h)}\\
 \le c h^2 \left(\norm{\omega \nabla v_h}_{L^2(D_{3\delta} \setminus D_\delta)} + \norm{v_h}_{L^2(D_{3\delta} \setminus D_\delta)} \right)
\end{multline*}
for every $v_h \in V_h$. The constant $c$ depends on $\delta$.
\end{lemma}
\begin{proof}
We have $\omega^2 v_h = v_h$ on $D_{2\delta}$ by definition of $\omega$. Since $h<\delta$ we have this equality on every cell $K \in \T_h$ with $K \cap D_\delta \neq \emptyset$. Thus, we obtain
\[
(\operatorname{id} - i_h)(\omega^2 v_h)|_K = 0 \quad \text{for every } K \text{ with } K \cap D_\delta \neq \emptyset.
\]
This results in
\[
\ltwonorm{(\operatorname{id} - i_h)(\omega^2 v_h)}^2 = \sum_{K \cap D_\delta = \emptyset} \ltwonormk{(\operatorname{id} - i_h)(\omega^2 v_h)}^2.
\]
Interpolation estimate on cell $K$ provides
\[
\begin{aligned}
\ltwonormk{(\operatorname{id} - i_h)(\omega^2 v_h)} &\le ch^2 \ltwonormk{\nabla^2(\omega^2 v_h)}\\
&\le ch^2 \left(\linfnormk{\nabla^2(\omega^2)} \ltwonormk{v_h} + \linfnorm{\nabla w} \ltwonormk{\omega \nabla v_h}\right)\\
& \le ch^2 \left(\delta^{-2} \ltwonormk{v_h} + \delta^{-1} \ltwonormk{\omega \nabla v_h}\right),
\end{aligned}
\]
since $\nabla^2 v_h = 0$ on $K$. To complete the proof we argue in the same way for the term $\ltwonorm{\nabla(\operatorname{id} - i_h)(\omega^2 v_h)}$ and use the fact that $\omega = 0$ on $\Omega \setminus D_{3\delta}$. Putting terms together results in the desired estimate.
\end{proof}

For the next results we again decompose the error $e = u-u_h$ as $e = \eta + \xi_h$ with
\begin{equation}\label{eq:eta_xih}
\eta = u - R_h u \quad \text{and} \quad \xi_h = R_h u - u_h \in V_h.
\end{equation}
We rewrite Galerkin orthogonality \cref{FEM:semilinear_Galerkin_orthogonality} using this decomposition resulting in
\[
	(\nabla \eta,\nabla \phi_h) +(\nabla \xi_h,\nabla \phi_h) + ( d(u)-d(u_{h}),\phi_h ) = 0 \quad \text{for all } \phi_h \in V_h
\]
and thus by definition of the Ritz projection
\begin{equation}\label{eq:gal_orth}
	(\nabla \xi_h,\nabla \phi_h) + ( d(u)-d(u_{h}),\phi_h ) = 0 \quad \text{for all } \phi_h \in V_h.
\end{equation}

\begin{lemma}\label{lemma:est_om_nab_xih}
	Let the nonlinearity $d$ fulfill the \cref{ass:1} and \cref{ass:2}. Let $\omega$ be defined in \cref{eq:omega} and let $\eta$ and $\xi_h \in V_h$ be given as in \cref{eq:eta_xih}. Then there holds
	\[
	\ltwonorm{\omega \nabla \xi_h}^2 \le ch^4 + \frac{1}{2\delta^2} \ltwonorm{\xi_h}^2 + c\norm{\xi_h}^2_{L^2(\Omega \setminus D_\delta)}.
	\] 
\end{lemma}
\begin{proof}
We obtain using Galerkin orthogonality \cref{eq:gal_orth}
\[
\begin{aligned}
\ltwonorm{\omega \nabla \xi_h}^2 &= (\nabla \xi_h, \omega^2 \nabla \xi_h) =(\nabla \xi_h, \nabla (\omega^2\xi_h)) - 2(\nabla \xi_h, \omega \xi_h \nabla \omega)\\
& = (\nabla \xi_h, \nabla i_h (\omega^2\xi_h)) +(\nabla \xi_h, \nabla(\operatorname{id} - i_h) (\omega^2\xi_h))- 2(\nabla \xi_h, \omega \xi_h \nabla \omega)\\
& = -( d(u)-d(u_{h}), i_h (\omega^2\xi_h))+(\nabla \xi_h, \nabla(\operatorname{id} - i_h) (\omega^2\xi_h))- 2(\nabla \xi_h, \omega \xi_h \nabla \omega)\\
&=-( d(u)-d(u_{h}), \omega^2\xi_h) + ( d(u)-d(u_{h}), (\operatorname{id} - i_h)(\omega^2\xi_h)) \\
&\qquad +(\nabla \xi_h, \nabla(\operatorname{id} - i_h) (\omega^2\xi_h))- 2(\nabla \xi_h, \omega \xi_h \nabla \omega).
\end{aligned}
\]
For the first term we get by the monotonicity of $d$
\[
\begin{aligned}
-( d(u)-d(u_{h}), \omega^2\xi_h) &= -( d(u)-d(R_h u), \omega^2\xi_h) - (d(R_h u) - d(u_h),\omega^2(R_h u - u_h))\\
& \le -( d(u)-d(R_h u), \omega^2\xi_h).
\end{aligned}
\]
This results in
\[
\begin{aligned}
\ltwonorm{\omega \nabla \xi_h}^2 &\le -( d(u)-d(R_h u), \omega^2\xi_h) + ( d(u)-d(u_{h}), (\operatorname{id} - i_h)(\omega^2\xi_h)) \\
&\qquad +(\nabla \xi_h, \nabla(\operatorname{id} - i_h) (\omega^2\xi_h))- 2(\nabla \xi_h, \omega \xi_h \nabla \omega)\\
&= J_1 +J_2 +J_3 +J_4.
\end{aligned}
\]
We estimate the four terms separately. From \cref{eq:D_4d} and the fact that
\[
\linfnorm{u-R_h u} \le ch^{2-\frac{N}{2}} \ltwonorm{\Lap u}
\]
we obtain
\[
\abs{R_h u(x)} \le \rho \quad \text{for all }\; x \in D_{4\delta}
\]
and thus by \cref{ass:2} we get
\[
\abs{d(x,u(x))-d(x,R_h u(x))} \le L \abs{\eta(x)} \quad \text{for almost all }\; x \in D_{4\delta}.
\]
Therefore, we obtain
\[
J_1 \le L \ltwonorm{\eta} \norm{\xi_h}_{L^2(D_{3\delta})} \le c h^4 + \frac{1}{2 \delta^2} \ltwonorm{\xi_h}^2,
\]
where we used estimate from \cref{lemma:standard_Rh} for $\eta$. For $J_2$ we use boundedness of $u$ and $u_h$ in $L^\infty(\Omega)$ and obtain
\[
\begin{aligned}
J_2 &\le c \ltwonorm{(\operatorname{id} - i_h)(\omega^2\xi_h)}
\le c h^2 \left(\norm{\omega \nabla \xi_h}_{L^2(D_{3\delta} \setminus D_\delta)} + \norm{\xi_h}_{L^2(D_{3\delta} \setminus D_\delta)} \right)\\
& \le ch^4 + \frac{1}{4} \ltwonorm{\omega \nabla \xi_h}^2 + c \norm{\xi_h}^2_{L^2(\Omega \setminus D_\delta)},
\end{aligned}
\]
where we used superapproximation estimate from \cref{lemma:superconv}. To estimate $J_3$ we observe that $(\operatorname{id} - i_h)(\omega^2\xi_h)$ vanishes on all cells $K$ with $K \cap D_\delta \neq \emptyset$, cf. proof of \cref{lemma:superconv}. Therefore we get by by inverse estimate and \cref{lemma:superconv}
\[
\begin{aligned}
J_3 &\le \norm{\nabla \xi_h}_{L^2(\Omega \setminus D_\delta)} \ltwonorm{\nabla(\operatorname{id} - i_h)(\omega^2\xi_h) }\\
& \le c h^{-1} \norm{\xi_h}_{L^2(\Omega \setminus D_\delta)} h \left(\norm{\omega \nabla \xi_h}_{L^2(D_{3\delta} \setminus D_\delta)} + \norm{\xi_h}_{L^2(D_{3\delta} \setminus D_\delta)} \right)\\
&\le \frac{1}{4} \ltwonorm{\omega \nabla \xi_h}^2  + c \norm{\xi_h}^2_{L^2(\Omega \setminus D_\delta)}.
\end{aligned}
\]
For $J_4$ we get
\[
J_4 \le \frac{1}{4} \ltwonorm{\omega \nabla \xi_h}^2 + c \norm{\xi_h}^2_{L^2(\Omega \setminus D_\delta)}.
\]
Putting terms together and absorbing  $\frac{3}{4} \ltwonorm{\omega \nabla \xi_h}^2$ into the left-hand side, we complete the proof.
\end{proof}

The following theorem provides a global estimate for $u-u_h$ with respect to the $L^2(\Omega)$ norm.
\begin{theorem}\label{theorem:L2}
Let the nonlinearity $d$ fulfill the \cref{ass:1} and \cref{ass:2}. Let $f \in L^\infty(\Omega)$ and $u \in H^1_0(\Omega)\cap C_0(\Omega)$ be the solution of
\cref{PDE:poisson_semilinear}. Let $u_h \in V_h$ be  the solution of
\cref{FEM:poisson_semilinear_weak_h}. Then, there holds
\[
\ltwonorm{u-u_h} \le c h^2 \lh^2.
\]
\end{theorem}
\begin{proof}
We use again the decomposition \cref{eq:eta_xih} for $e=u-u_h = \eta+\xi_h$. For $\eta$ we have by \cref{lemma:standard_Rh}
\[
\ltwonorm{\eta} \le ch^2 \ltwonorm{\Lap u}.
\]
To estimate $\xi_h$ we use the Poincare inequality on $D_\delta$
\[
\norm{v}_{L^2(D_\delta)} \le \delta \norm{\nabla v}_{L^2(D_\delta)} \quad \text{for every } v \in H^1_0(\Omega).
\]
Thus, we have
\[
\begin{aligned}
\ltwonorm{\xi_h}^2 &= \norm{\xi_h}^2_{L^2(D_\delta)} + \norm{\xi_h}^2_{L^2(\Omega \setminus D_\delta)}\\
& \le \delta^2 \norm{\nabla \xi_h}^2_{L^2(D_\delta)}  + c\norm{\xi_h}^2_{L^2(\Omega \setminus D_\delta)}\\
& \le \delta^2 \ltwonorm{\omega \nabla \xi_h}^2 + c\norm{\xi_h}^2_{L^2(\Omega \setminus D_\delta)}.
\end{aligned}
\]
Using \cref{lemma:est_om_nab_xih} and \cref{theorem:interior} we get
\[
\begin{aligned}
\ltwonorm{\xi_h}^2 &\le ch^4 + \frac{1}{2} \ltwonorm{\xi_h}^2 + c\norm{\xi_h}^2_{L^2(\Omega \setminus D_\delta)}\\
 \le & ch^4 + \frac{1}{2} \ltwonorm{\xi_h}^2 + c\norm{\xi_h}^2_{L^\infty(\Omega \setminus D_\delta)}\\
 &\le ch^4 \lh^4 + \frac{1}{2} \ltwonorm{\xi_h}^2.
\end{aligned}
\] 
Absorbing last term into the left-hand side completes the proof.
\end{proof}

\begin{remark}
If the boundary of $\Omega$ is assumed to be smooth, then the above result can be shown without using \cref{ass:2}. The question whenever it is possible to show this (or similar) estimate on a convex polygonal/polyhedral domain without \cref{ass:2} seems to be open.
\end{remark}

\section{Numerical examples}
In this section we present a numerical example illustration our results. We choose a convex polygonal domain $\Omega \subset \R^2$ with the largest angle $\omega_\Omega = \nicefrac{3\pi}{4}$. The nonlinearity is set as
\[
d(x,u) = r(u+1), \quad r(t) = 50 \sgn(t) \abs{t}^\frac{1}{3}
\]
and the right-hand side $f=1$. The expected converges results are
\[
\ltwonorm{u-u_h} = {\mathcal O}(h^2\lh^2) \quad\text{and}\quad \linfnorm{u-u_h} = {\mathcal O}(h^\gamma\lh), \; \gamma < \frac{4}{3}.
\]
The following plot presents the results.
\begin{figure}
\begin{tikzpicture}
	\begin{loglogaxis} [%
		xlabel={mesh size $h$},
		grid=major,
		legend pos=south east,
		legend cell align=left,
		width=\textwidth,
		height=0.6\textwidth,
		cycle list name=black white,
		]
		
		\addplot table [x=h, y=el2, col sep=semicolon] {NumEx/c.csv};
		\addlegendentry{$\ltwonorm{u-u_h}$}
		\addplot table [x=h, y=el8, col sep=semicolon] {NumEx/c.csv};
		\addlegendentry{$\linfnorm{u-u_h}$}
		
		\logLogSlopeTriangle{0.25}{0.15}{0.36}{1.333333333333}{black}{h^{\frac43}}
		\logLogSlopeTriangle{0.25}{0.15}{0.073}{2}{black}{h^2}
	\end{loglogaxis}
\end{tikzpicture}
\caption{Error behavior for $\ltwonorm{u-u_h}$ and $\linfnorm{u-u_h}$ under mesh refinement}
\end{figure}
\section*{Acknowledgments}
We would like to thank Dr. Dominik Meidner for scientific exchange and for helping
us to prepare the numerical example.

\bibliographystyle{siamplain}
\bibliography{lit}

\begin{thebibliography}{10}

\bibitem{Barrett:1992}
{\sc J.~W. Barrett}, {\em Finite element approximation of a non-{L}ipschitz
  nonlinear eigenvalue problem}, RAIRO Mod\'el. Math. Anal. Num\'er., 26
  (1992), pp.~627--656.

\bibitem{BarrettShanahan:1991}
{\sc J.~W. Barrett and R.~M. Shanahan}, {\em Finite element approximation of a
  model reaction-diffusion problem with a non-{L}ipschitz nonlinearity}, Numer.
  Math., 59 (1991), pp.~217--242.

\bibitem{BrennerScott:2008}
{\sc S.~C. Brenner and L.~R. Scott}, {\em The Mathematical Theory of Finite
  Element Methods}, vol.~15 of Texts in Applied Mathematics, Springer, New
  York, third~ed., 2008.

\bibitem{Casas:1993}
{\sc E.~Casas}, {\em Boundary control of semilinear elliptic equations with
  pointwise state constraints}, SIAM J. Control Optim., 31 (1993),
  pp.~993--1006.

\bibitem{CasasMateos:2002}
{\sc E.~Casas and M.~Mateos}, {\em Uniform convergence of the {FEM}.
  {A}pplications to state constrained control problems}, vol.~21, 2002,
  pp.~67--100.
\newblock Special issue in memory of Jacques-Louis Lions.

\bibitem{Christof:2024}
{\sc C.~Christof}, {\em Optimal control of semilinear elliptic partial
  differential equations with non-lipschitzian nonlinearities}.
\newblock {a}rXiv:2406.03110 [math.NA], 2024,
  \url{https://doi.org/10.48550/arXiv.2406.03110},
  \url{https://doi.org/10.48550/arXiv.2406.03110}.

\bibitem{FrehseRannacher:1976}
{\sc J.~Frehse and R.~Rannacher}, {\em Eine {$L^1$}-{F}ehlerabsch\"{a}tzung
  f\"{u}r diskrete {G}rundl\"{o}sungen in der {M}ethode der finiten
  {E}lemente}, in Finite {E}lemente ({T}agung, {I}nst. {A}ngew. {M}ath.,
  {U}niv. {B}onn, {B}onn, 1975), vol.~89 of Bonner Math. Schriften,
  Universit\"{a}t Bonn, Institut f\"{u}r Angewandte Mathematik, Bonn, 1976,
  pp.~92--114.

\bibitem{GilbargTrudinger:2001}
{\sc D.~Gilbarg and N.~S. Trudinger}, {\em Elliptic Partial Differential
  Equations of Second Order}, Classics in Mathematics, Springer-Verlag, Berlin,
  2001.
\newblock Reprint of the 1998 edition.

\bibitem{Grisvard:1985}
{\sc P.~Grisvard}, {\em Elliptic Problems in Nonsmooth Domains}, vol.~24 of
  Monographs and Studies in Mathematics, Pitman (Advanced Publishing Program),
  Boston, MA, 1985.

\bibitem{GuzmanLeykekhmanannRossmannSchatz:2009}
{\sc J.~Guzm\'{a}n, D.~Leykekhman, J.~Rossmann, and A.~H. Schatz}, {\em
  H\"{o}lder estimates for {G}reen's functions on convex polyhedral domains and
  their applications to finite element methods}, Numer. Math., 112 (2009),
  pp.~221--243.

\bibitem{HafemeyerKahlePfefferer:2020}
{\sc D.~Hafemeyer, C.~Kahle, and J.~Pfefferer}, {\em Finite element error
  estimates in {$L^2$} for regularized discrete approximations to the obstacle
  problem}, Numer. Math., 144 (2020), pp.~133--156.

\bibitem{KinderlehrerStampacchia:2000}
{\sc D.~Kinderlehrer and G.~Stampacchia}, {\em An Introduction to Variational
  Inequalities and their Applications}, vol.~31 of Classics in Applied
  Mathematics, Society for Industrial and Applied Mathematics (SIAM),
  Philadelphia, PA, 2000.
\newblock Reprint of the 1980 original.

\bibitem{KnabnerRannacher:2017}
{\sc P.~Knabner and R.~Rannacher}, {\em A priori error analysis for the
  {G}alerkin finite element semi-discretization of a parabolic system with
  non-{L}ipschitzian nonlinearity}, Vietnam J. Math., 45 (2017), pp.~179--198.

\bibitem{LeykekhmanD_VexlerB_2016c}
{\sc D.~Leykekhman and B.~Vexler}, {\em Finite element pointwise results on
  convex polyhedral domains}, SIAM J. Numer. Anal., 54 (2016), pp.~561--587.

\bibitem{NitscheSchatz:1974}
{\sc J.~A. Nitsche and A.~H. Schatz}, {\em Interior estimates for
  {R}itz-{G}alerkin methods}, Math. Comp., 28 (1974), pp.~937--958,
  \url{https://doi.org/10.2307/2005356}, \url{https://doi.org/10.2307/2005356}.

\bibitem{Nochetto:1988}
{\sc R.~H. Nochetto}, {\em Sharp {$L^\infty$}-error estimates for semilinear
  elliptic problems with free boundaries}, Numer. Math., 54 (1988),
  pp.~243--255.

\bibitem{Ruzicka:2020}
{\sc M.~R\r{u}\v{z}i\v{c}ka}, {\em Nichtlineare Funktionalanalysis},
  Masterclass, Springer Spektrum Berlin, Heidelberg, second~ed., 2020.

\bibitem{SchatzWahlbin:1977}
{\sc A.~H. Schatz and L.~B. Wahlbin}, {\em Interior maximum norm estimates for
  finite element methods}, Math. Comp., 31 (1977), pp.~414--442.

\bibitem{Troeltzsch:2010}
{\sc F.~Tr\"{o}ltzsch}, {\em Optimal control of partial differential
  equations}, vol.~112 of Graduate Studies in Mathematics, American
  Mathematical Society, Providence, RI, 2010.
\newblock Theory, methods and applications, Translated from the 2005 German
  original by J\"{u}rgen Sprekels.

\bibitem{Zeidler:IIB:1990}
{\sc E.~Zeidler}, {\em Nonlinear Functional Analysis and its Applications.
  II/B}, Springer-Verlag, New York, 1990.
\newblock Nonlinear monotone operators, Translated from the German by the
  author and Leo F. Boron.

\end{thebibliography}
\end{document}